\documentclass[11pt]{article}
\title{On the lower bound of the spectral norm of symmetric random  matrices with independent entries}
\author{Sandrine P\'{e}ch\'{e} \thanks{Institut Fourier BP 74, 100 Rue des maths, 38402 Saint Martin d'Heres, France  
(permanent address) and
Department of Mathematics,
University of California at Davis, 
One Shields Ave., Davis, CA 95616, USA.
E-mail: Sandrine.Peche@ujf-grenoble.fr}
\and Alexander Soshnikov \thanks{
Department of Mathematics,
University of California at Davis, 
One Shields Ave., Davis, CA 95616, USA.
E-mail address: soshniko@math.ucdavis.edu.
Research was supported in part by the 
NSF grants DMS-0405864 and DMS-0707145. 
}  }
\usepackage{epsfig}
\usepackage{amsfonts}
\usepackage{amssymb}
\usepackage{amsthm}
\usepackage{amstext}
\usepackage{amscd}
\usepackage{amsmath}
\usepackage{delarray}
\begin{document}
\maketitle
\newtheorem{theo}{Theorem}[section]
\newtheorem{prop}{Proposition}[section]
\newtheorem{lemme}{Lemma}[section]
\newtheorem{conjecture}{Conjecture}[section]
\newtheorem{definition}{Definition}[section]
\newtheorem{fact}{Fact}[section]
\newtheorem{hyp}{Assumption}[section]
\theoremstyle{remark}
\newtheorem{rem}{Remark}
\newtheorem{remark}{Remark}[section]
\newtheorem{Remark}{Remark}[section]
\newtheorem{Notationnal remark}{Remark}[section]
\newcommand{\bremnot}{\begin{Notationnal remark}}
\newcommand{\eremnot}{\end{Notationnal remark}}
\newcommand{\brem}{\begin{remark}}
\newcommand{\erem}{\end{remark}}
\newcommand{\bconj}{\begin{conjecture}}
\newcommand{\econj}{\end{conjecture}}
\newcommand{\bdefi}{\begin{definition}}
\newcommand{\edefi}{\end{definition}}
\newcommand{\bt}{\begin{theo}}
\newcommand{\bfa}{\begin{fact}}
\newcommand{\efa}{\end{fact}}
\newcommand{\Si}{\Sigma}
\newcommand{\mbE}{\mathbb{E}}
\newcommand{\mL}{\mathcal{L}}
\newcommand{\mP}{\mathcal{P}}
\newcommand{\mQ}{\mathcal{Q}}
\newcommand{\mR}{\mathcal{R}}
\newcommand{\et}{\end{theo}}
\newcommand{\bp}{\begin{prop}}
\newcommand{\ep}{\end{prop}}
\newcommand{\bl}{\begin{lemme}}
\newcommand{\el}{\end{lemme}}
\newcommand{\be}{\begin{equation}}
\newcommand{\ee}{\end{equation}}
\newcommand{\bmp}{\begin{pmatrix}}
\newcommand{\emp}{\end{pmatrix}}
\newcolumntype{L}{>{$}l<{$}}
\newenvironment{Cases}{\begin{array}\{{lL.}}{\end{array}}
\begin{abstract}
We show that the spectral radius of an $N\times N$ random symmetric matrix with i.i.d. bounded centered but non-symmetrically
distributed entries  is bounded from below by \\
$ 2 \*\sigma - o( N^{-6/11+\varepsilon}), $ where $\sigma^2 $ is the variance of the matrix 
entries and $\varepsilon $ is an arbitrary small positive number.   Combining with our previous result from [7], this proves 
that for any  $\varepsilon >0, \ $ one has
$ \|A_N\| =2 \*\sigma + o( N^{-6/11+\varepsilon}) $ with probability going to $ 1 $ as $N \to \infty. $  
\end{abstract}

\section{Introduction}

Wigner random matrices were introduced by E.Wigner about fifty years ago (\cite{Wig1}, \cite{Wig2}) as a model to study the 
statistics of resonance levels for neutrons off heavy nuclei.  Nowadays, there are many fruitful connections between Random 
Matrix Theory and Mathematical Physics, Probability Theory, Integrable Systems, Number Theory, Quantum Chaos, Theoretical 
Computer Theory, Combinatorics, Statistics, and many other areas of science.

Let $A_N$ be a sequence of real symmetric Wigner random matrices with non symmetrically distributed entries.  In other words,
$$A_N=\frac{1}{\sqrt N}\left ( a_{ij}\right)_{i,j=1}^N,$$ where the $a_{ij}, i\leq j$ are i.i.d. random variables such that
\be \label{H3} 
\mbE a_{ij}=0, \, \mbE a_{ij}^2=\sigma^2, \,  \mbE a_{ij}^3=\mu_3, \text{ and } |a_{ij}| \leq C, \  
\forall 1 \leq i,j \leq N,
\ee
where $C$ is some positive constant that does not depend on $N.$  The common third moment $\mu_3$ is not necessarily zero, 
which allows us to study the case when the marginal distribution of matrix entries is not symmetric.

Let us denote by $\|A_N\| $ the spectral norm of the matrix $A_N, \ 
\|A_N\|=max_{1\leq i \leq N} |\lambda_i|, $ where $\lambda_1, \ldots,
\lambda_N$ are the eigenvalues of $A_N. \ $  Clearly, the eigenvalues of $A_N$ are real random variables.
It was proved in \cite{PecheSos} that for an arbitrary small positive number $\varepsilon >0$
the spectral norm of  $A_N$ is bounded as
\begin{equation}
\label{ocenka17}
\|A_N \| \leq 2 \*\sigma + o( N^{-6/11+\varepsilon}), 
\end{equation}
with probability going to 1.
In this paper, we prove that $ 2 \*\sigma + o( N^{-6/11+\varepsilon}) $ is also a lower bound for $ \|A_N\|. $
The main result of the paper is the following

\bt \label{mt}
Let $ \|A_N\|$ denote the spectral norm of the matrix $A_N$ and $ \varepsilon >0$.  Then
\begin{equation}
\label{ocenka100}
 \|A_N\| \geq 2 \*\sigma - o(N^{-6/11+\varepsilon}) 
\end{equation}
with probability going to 1 as $ N \to \infty. $
\et

Combining the result of Theorem 1.1 with  (\ref{ocenka17}), we obtain

\bt \label{mt2}
Let $ \|A_N\|$ denote the spectral norm of the matrix $A_N$ and $ \varepsilon >0$.  Then
\begin{equation}
\label{ocenka}
 \|A_N\| = 2 \*\sigma + o(N^{-6/11+\varepsilon}) 
\end{equation}
with probability going to 1 as $ N \to \infty. $
\et

\brem \label{rem: indepen}  In fact, one does not need the assumption that the matrix entries are identically distributed as 
long as $ \{ a_{ij}, 1\leq i \leq j \leq N \} $ are independent, uniformly bounded centralized random variables with the same 
variance $\sigma^2 $ off the diagonal.  The proofs of the results of the present paper and of \cite{PecheSos} still hold 
without any significant alterations since we only use the upper bounds $|\mbE a_{ij}^k| \leq C^k $ on the third and higher 
moments, i.e. for $k\geq 3, $ and not the exact values of these moments.
\erem

\brem \label{rem: hermitian} Similar results hold for Hermitian Wigner matrices as well.  Since the proof is essentially 
the same,  we will discuss only the real symmetric case in this paper.
\erem
We remark that $2 \*\sigma $ is the right edge of the support of the Wigner semicircle law, and, therefore, it immediately 
follows from the classical result of Wigner (\cite{Wig1}, \cite{Wig2}, \cite{Arn}) that
for any fixed $\delta >0, \  \mathbb{P} \left(  \|A_N\| \geq 2 \*\sigma -\delta \right) \to 1 $ as $ N \to \infty. $

A standard way to obtain an upper bound on the spectral norm is to study the asymptotics of $ \mbE [ Tr A_N^{2s_N} ] $ for 
integers $ s_N $ proportional to $N^{\gamma}, \ \gamma >0. $  If one can show that
\begin{equation}
\label{ocenka1}
\mbE [ Tr A_N^{2s_N} ] \leq Const_1 \* N^{\gamma_1}\*(2\sigma)^{2s_N}, 
\end{equation}
where $ s_N= Const \*N^{\gamma}\*(1+o(1)), $ and $ Const_1 $ and $\gamma_1$ depend only on $Const $ and $\gamma, $ 
one can prove that
\begin{equation}
\label{ocenka2}
\|A_N\| \leq 2 \*\sigma + O( N^{-\gamma} \*\log N) 
\end{equation}
with probability going to 1 by
using the upper bound $ \mbE [ \|A_N\|^{2s_N}  ] \leq \mbE [ Tr A_N^{2s_N} ] $ and the Markov inequality.   In particular,  
F\"{u}redi and Koml\'{o}s in \cite{Fur-Kom} were able to prove (\ref{ocenka2}) for $\gamma \leq 1/6, $ and Vu  \cite{Vu}
extended their result to  $\gamma \leq 1/4. $  Both  papers \cite{Fur-Kom} and \cite{Vu} treated the case when the matrix 
entries $\{a_{ij}\} $ are uniformly bounded.
In  \cite{PecheSos}, we were able to prove that
\begin{equation}
\label{ocenka101}
\mbE [ Tr A_N^{2s_N} ] =\frac{N}{\pi^{1/2}\*s_N^{3/2}}\*(2\sigma)^{2s_N}\* (1+o(1)), 
\end{equation}
for $s_N =O(N^{6/11 -\varepsilon}) $ and any $\varepsilon >0, $ thus establishing (\ref{ocenka17}).  Again, we restricted our 
attention in  \cite{PecheSos} to the case of uniformly bounded entries.
The proof relies on combinatorial arguments going back to \cite{Si-So1}, \cite{Si-So2}, and \cite{So1}.

More is known if the matrix entries of
a Wigner matrix have symmetric distribution (so, in particular, the odd moments of matrix entries vanish).
In the case of symmetric marginal distribution of matrix entries, one can relax the condition that
$ (a_{ij}) $ are uniformly bounded and assume that the marginal distribution is sub-Gaussian.
It was shown by Tracy and Widom in \cite{TW1} in the Gaussian (GOE) case that
the largest eigenvalue deviates from the soft edge $2\sigma $
on the order $ O(N^{-2/3}) $ and the limiting distribution of the rescaled largest eigenvalue  obeys
Tracy-Widom law (\cite{TW1}):
$$  \lim_{N\to \infty} \mathbb{P} \left( \lambda_{max} \leq 2\sigma + \sigma \* x \*N^{-2/3} \right)= 
\exp\left(-1/2\*\int_x^{\infty} q(t)+(t-x)\*q^2(t)\*dt \right), $$  
where $q(x) $ is the solution of the Painl\'{e}ve II differential equation $ \ q''(x)=x\*q(x) +2\*q^3(x) \ $ with the 
asymptotics  at infinity
$ q(x) \sim Ai(x) $ as $ x \to +\infty. \ $   It was shown in \cite{So1} that this behavior is universal for Wigner matrices 
with sub-Gaussian and symmetrically distributed entries.  Similar results hold in the Hermitian case 
(see \cite{TW2}, \cite{So1}).
It is reasonable to expect that in the non-symmetric case, the largest eigenvalue 
will have the Tracy-Widom distribution in the limit as well.

The lower bonds on the spectral norm of a Wigner random matrix with non-symmetrically distributed entries were probably 
considered to be more difficult than the upper bounds.  Let us again restrict our attention
to the case when matrix entries are uniformly bounded.
It was claimed in  \cite{Fur-Kom} that the estimate of the type (\ref{ocenka1}) for
$\gamma \leq 1/6 $ immediately implies the lower bound $$\|A_N\| \geq 2 \*\sigma - O(N^{-1/6}\*\log N). $$  As noted by 
Van Vu in \cite{Vu}, ``We do not see any way to materialize this idea.''  We concur with this opinion.  In the next section, 
we show that (\ref{ocenka1}) implies a rather weak estimate
\begin{equation}
\label{ocenka3}
\mathbb{P} \left(\|A_N\| \geq 2 \*\sigma -N^{-6/11 +\delta} \right) \geq N^{-9/11 +\delta}, 
\end{equation}
for small $\delta>0 $ and sufficiently large $N.$   Combining (\ref{ocenka3}) with the concentration 
of measure inequalities for $\|A_N\| $ (see \cite{GZ}, \cite{AKV}),  one then obtains that for Wigner matrices with uniformly 
bounded entries
\begin{equation}
\label{non-optimal}
\mathbb{P} \left(\|A_N\| \geq 2 \*\sigma \*(1 - C \* N^{-1/2}\*\sqrt{\log N}) \right) \to 1 \text{ as } N \to \infty,
\end{equation}
where $C$ is the same as in (\ref{H3}).
The proof of Theorem 1.1 will
be given in Section 3, where we establish the analogue of the law of large numbers for 
$ Tr A_N^{2s_N}$ with $ s_N=O(N^{6/11-\varepsilon}), $ proving that 
\begin{equation}
\label{LLN}
Tr A_N^{2s_N}  =\frac{N}{\pi^{1/2}\*s_N^{3/2}}\*(2\sigma)^{2s_N}\* (1+o(1))
\end{equation}
with probability going to 1 as $N \to \infty. $

\section{Preliminary lower bound}

Without loss of generality, we can assume $\sigma=1/2. $  This conveniently sets the right edge of the Wigner semicircle law
to be $1.$  Let us fix $ 0< \delta < 6/11, $ and denote $ \Omega_N= \{ \|A_N\|>1 -N^{-6/11+\delta} \}. $  Choose $s_N$ to
be an integer such that 
$s_N=N^{6/11 -\varepsilon} \*(1+o(1)) $ and $2 \*\delta/3 < \varepsilon <\delta. $  Let us denote by $1_{\Omega} $ the 
indicator of 
the set $\Omega $  and by $\Omega^c$ the complement of $\Omega. $ Then
\begin{equation}
\label{brest}
\mbE \left( Tr A_N^{2s_N}\*1_{\Omega_N^c} \right) \leq N \* (1-N^{-6/11+\delta})^{2s_N} 
\leq  N \* \left(\exp(-N^{-6/11+\delta}) \right)^{N^{6/11 -\varepsilon}}= O\left(N \* e^{-N^{\delta-\varepsilon}}\right)
\end{equation}
which is $o(1)$ as $N\to \infty.$  Let us now partition $\Omega_N$ as the disjoint union $\Omega_N=\Omega_N^1 \bigsqcup
\Omega_N^2 , $ where
$$ \Omega_N^1= \{ 1 -N^{-6/11+\delta} < \|A_N\| < 1 +N^{-6/11+\varepsilon} \} \text{ and }
\Omega_N^2= \{ \|A_N\| \geq 1 +N^{-6/11+\varepsilon} \}. $$
Then
\begin{equation}
\label{minsk}
\mbE \left( Tr A_N^{2s_N}\*1_{\Omega_N^1} \right) \leq N \* (1+N^{-6/11+\varepsilon})^{2s_N}  \* 
\mathbb{P} \left(\Omega_N^1\right)
\leq  N\* (e^2+o(1))\*  \mathbb{P} \left(\Omega_N^1\right).
\end{equation}
As for $\mbE \left( Tr A_N^{2s_N}\*1_{\Omega_N^2} \right), $ one can show that
\begin{eqnarray}
\label{smolensk}
& & \mbE \left( Tr A_N^{2s_N}\*1_{\Omega_N^2} \right) \leq \mbE \left(N \* \|A_N\|^{2 s_N} \* 1_{\Omega_N^2}\right) \leq 
\mbE \left( N \* \|A_N\|^{2\*[ N^{6/11 -\varepsilon/4}]} \* 1_{\Omega_N^2} \right) \leq \\
& & N \* (1 +N^{-6/11 +\varepsilon})^{- 2\* [N^{6/11-\varepsilon/4 }]} \* 
\mbE \left( \|A_N\|^{4 \* [N^{6/11-\varepsilon/4}]} \* 1_{\Omega_N^2} \right) \leq N \* e^{-2\*N^{3\*\varepsilon/4}} \*
\mbE \left( \|A_N\|^{4 \* [N^{6/11-\varepsilon/4}]} \right) \nonumber \\
& & \leq N^2\* e^{-2\*N^{3\*\varepsilon/4}},  \nonumber
\end{eqnarray}
where in the last inequality we used (\ref{ocenka101}) (for $\sigma=1/2$) to get
$$ \mbE \left( \|A_N\|^{4 \* [N^{6/11-\varepsilon/4}]} \right) \leq  \mbE \left(Tr A_N^{4 \* [N^{6/11-\varepsilon/4}]} \right)
=\frac{N}{\pi^{1/2}\* (4\* N^{6/11-\varepsilon/4})^{3/2}}\*(1+o(1))\leq N
$$
for large $N.$

Combining the above estimates and (\ref{ocenka101}) (for $\sigma=1/2$), we obtain for sufficiently large $N$ that
\begin{equation}
\label{vitebsk}
\frac{N}{2 \*s_N^{3/2}} \leq 
\mbE \left( Tr A_N^{2s_N}\right) \leq O\left(N \* e^{-N^{\delta-\varepsilon}}\right) + O( N^2\* e^{-N^{3\*\varepsilon/4}})
+ \mathbb{P} \left(\Omega_N^1\right) \*  N\* (e^2+o(1)).
\end{equation}
Therefore,
\begin{equation}
\label{lida}
\mathbb{P} \left( \|A_N\|>1 -N^{-6/11+\delta} \right) \geq \mathbb{P} \left(\Omega_N^1\right) \geq \frac{e^{-2}}{2} 
\*s_N^{-3/2} \* (1+o(1))=N^{-9/11 + 3\*\varepsilon/2} \*(e^{-2}/2+o(1)) \geq N^{-9/11 +\delta}
\end{equation}
for sufficiently large $N$ (depending on $\delta$.)

It was shown by  by Alon, Krivelevich, and Vu (\cite{AKV}), and Guionnet and Zeitouni (\cite{GZ}) that for Wigner random 
matrices with bounded entries, the spectral norm is strongly concentrated around its mean.  Indeed, the spectral norm is a 
1-Lipschitz function of the matrix entries since 
$$ | \|A\|-\|B\| |\leq \|A-B\| \leq \|A-B\|_{HS} = \left(Tr \left( (A-B)\*(A-B)^t \right) \right)^{1/2}=
\left( \sum_{ij} (a_{ij}-b_{ij})^2 \right)^{1/2}, $$
where $\| \|_{HS} $ denotes the Hilbert-Schmidt norm.
Therefore, one can apply the concentration of measure results (\cite{Tal}, \cite{Led1}, 
\cite{Led}).
In particular (see Theorem 1 in (\cite{AKV})),
\begin{equation}
\label{grodno}
\mathbb{P} \left( | \|A_N\|- \mbE \|A_N\| |> C\* t \* N^{-1/2} \right) \leq 4 \* e^{-t^2/32},
\end{equation}
uniformly in $N$ and $t, $
where the constant $C$ is the same as in (\ref{H3}).  Combining (\ref{lida}) and (\ref{grodno}), we arrive at 
\begin{equation}
\label{redkino100}
 \mbE \|A_N\| \geq 1- \frac{3}{\sqrt{11}} \*C \*N^{-1/2}\*\sqrt{\log N}
\end{equation}
for sufficiently large $N$.  The last inequality together with (\ref{grodno}) then implies
(\ref{non-optimal}) (recall that we set $\sigma=1/2.$)

\section{Law of Large Numbers}

The main technical result of this section is the following analogue of the Law of Large Numbers for $Tr A_N^{2s_N}.$
\bp \label{prop: lln}
Let $s_N = O(N^{6/11 - \varepsilon}), $ where $ \varepsilon $ is an arbitrary small constant.  Then
\begin{equation}
\label{polock}
Tr A_N^{2s_N} = \mbE \left(Tr A_N^{2s_N} \right)    \* (1+ \delta_N), 
\end{equation}
where $ \mathbb{P} \left( |\delta_N| \geq N^{-1/22} \right) \to 0 $ as $ N \to \infty. $
\ep

The Proposition 3.1 combined with (\ref{ocenka101}) immediately implies that
\begin{equation}
\label{mai}
Tr A_N^{2s_N} = \frac{N}{\pi^{1/2}\*s_N^{3/2}}\*(2\sigma)^{2s_N}\* (1+ o(1)),
\end{equation}
with probability going to 1 as $N \to \infty.$  To make (\ref{mai}) more precise,  we can say that
the ratio of the l.h.s. and the r.h.s. of 
(\ref{mai}) goes to 1 in probability as $N \to \infty.$

The main part of the proof of Proposition  \ref{prop: lln} is the following bound on the variance.
\bl \label{lem: var}
Let $s_N = O(N^{6/11 - \varepsilon}), $ where $ \varepsilon $ is an arbitrary small constant.  There there exists $Const >0 $ 
such that
\begin{equation}
\label{pinsk}
\text{Var } (Tr A_N^{2s_N}) \leq Const \* \sqrt{s_N}\* (2\sigma)^{4s_N}.
\end{equation}
\el
Lemma is proven in the subsection below.  
Assuming Lemma \ref{lem: var}, we  obtain the proof of Proposition \ref{prop: lln} via the Chebyshev inequality.  Indeed,
it follows from (\ref{pinsk}) and (\ref{ocenka101}) that
\begin{equation}
\label{msbmsb}
\text{Var } \left(\frac{Tr A_N^{2s_N}}{\mathbb{E} [Tr A_N^{2s_N}]} \right) \leq \frac{Const \* \sqrt{s_N}\* (2\sigma)^{4s_N}}
{N^2 \* (2\sigma)^{4s_N} / (\pi \* s_N^3)} \* (1+o(1)) = O(s_N^{7/2} \*N^{-2})= 
O(N^{-1/11 -7\*\varepsilon/2}).
\end{equation}  
To finish the proof of the 
main result of the paper, we fix an arbitrary small positive constant $\delta>0$ and choose another constant $\varepsilon$ in 
such a way that $0 < \varepsilon < \delta. \ $  Setting $\sigma=1/2,$ we scale the eigenvalues in such a way that
the right edge of the Wigner semicircle law is equal to $1.$
Let us denote, as before, $ \Omega_N= \{ \|A_N\|>1 -N^{-6/11+\delta} \}. $
Choosing $s_N=N^{6/11 -\varepsilon} \*(1+o(1)), $  we note that on  $ \Omega_N^c $ 
$$ Tr A_N^{2s_N}\*1_{\Omega_N^c} \leq N \* (1-N^{-6/11+\delta})^{2s_N} 
= O\left(N \* e^{-N^{\delta-\varepsilon}}\right)=o(1). $$
At the same time, Proposition  \ref{prop: lln} implies (see (\ref{mai})) that
$ Tr A_N^{2s_N}  \geq N^{2/11} \ $ with probability going to 1.  Therefore,
\begin{equation}
\label{mogilev1}
\mathbb{P} \left( \|A_N\| \leq 1 -N^{-6/11+\delta} \right) \to 0  \text{ as } N \to \infty.
\end{equation}

\subsection{The proof of Lemma}

We now turn our attention to the variance of the trace, which can be considered as follows.  To express 
$\text{Var }\text{Tr}A_N^{2s_N} $ in 
terms of the matrix entries, we first write $\text{Tr}A_N^{2s_N}$ as the sum of the products of matrix entries, namely
we express $\text{Tr}A_N^{2s_N}$ as the sum of the diagonal entries of the matrix $A_N^{2s_N}.\ $  Therefore,
\begin{equation}
\label{sledy100}
\mathbb{E} \text{Tr}A_N^{2s_N}= \sum_{1\leq i_0, \ldots, i_{2\*s_N-1} \leq N} \mathbb{E} \prod_{0\leq k\leq 2\*s_N -1} 
a_{i_k i_{k+1}},
\end{equation}
where we assume that $i_{2s_N}=i_0.$  We can then rewrite $\mathbb{E} \text{Tr}A_N^{2s_N} $ as the sum over the set of closed
paths $\mathcal P = \{i_0 \to i_1 \to \ldots i_{2s_N-1}\to i_0 \} $ on the complete graph on the $N$ vertices 
$\{1,2,\ldots, N\}$ as
\begin{equation}
\label{sledy}
\mathbb{E} \text{Tr}A_N^{2s_N}= \sum_{ \mathcal P} \mathbb{E} \prod_{(i_k i_{k+1}) \in \mathcal P} 
a_{i_k i_{k+1}}.
\end{equation}

In a similar fashion (again using the agreement that $i_{2\*s_N}=i_0 $ and $ j_{2\*s_N}=j_0 $ ), we can write
\begin{eqnarray}
&&\text{Var }\text{Tr}A_N^{2s_N}= 
\frac{1}{N^{2s_N}} \sum_{1\leq i_0, \ldots, i_{2\*s_N-1} \leq N} \ \sum_{1\leq j_0, \ldots, j_{2\*s_N-1} \leq N} \crcr
&&\Big [\mathbb{E} \prod_{0\leq k\leq 2\*s_N -1} \ \prod_{0\leq l\leq 2\*s_N -1}  a_{i_k i_{k+1}} a_{j_l j_{l+1}}  -
\mathbb{E}\prod_{0\leq k\leq 2\*s_N -1}  a_{i_k i_{k+1}}
\* \mathbb{E} \prod_{0\leq l\leq 2\*s_N -1}
a_{j_l j_{l+1}} \Big ] \nonumber \crcr
&&=\frac{1}{N^{2s_N}}\sum_{\mathcal P_1,\mathcal P_2}\Big [\mathbb{E}\prod_{(i_ki_{k+1})\in \mathcal P_1}
\prod_{(j_l j_{l+1})\in \mathcal P_2} a_{i_k i_{k+1}}a_{j_l j_{l+1}}-\mathbb{E}\prod_{(i_ki_{k+1})\in \mathcal P_1}
a_{i_ki_{k+1}}\mathbb{E}\prod_{(j_lj_{l+1})\in \mathcal P_2} a_{j_l j_{l+1}}\Big ] \nonumber \crcr
\label{belovezh}
&&=\frac{1}{N^{2s_N}}\sum_{\mathcal P_1,\mathcal P_2}^{\star}\Big [\mathbb{E}\prod_{(i_ki_{k+1})\in \mathcal P_1}
\prod_{(j_l j_{l+1})\in \mathcal P_2} a_{i_ki_{k+1}}a_{j_l j_{l+1}} - \mathbb{E}\prod_{(i_ki_{k+1})\in \mathcal P_1}
a_{i_ki_{k+1}}\mathbb{E}\prod_{(j_l j_{l+1})\in \mathcal P_2} a_{j_l j_{l+1}} 
\Big], 
\end{eqnarray}
where $\mathcal P_1 $ and $\mathcal P_2$ are closed paths of length $2s_N, $
$$ \mathcal P_1= \{i_0 \to i_1 \to \ldots i_{2s_N-1}\to i_0 \}  \ \ \ \text{and} \ \ \ 
 \mathcal P_2= \{j_0 \to j_1 \to \ldots j_{2s_N-1}\to j_0 \}. $$
The starred summation symbol $\sum_{\mathcal P_1,\mathcal P_2}^{\star} $ in the last line of the previous array of equations
means that the summation is restricted to the set of 
the pairs of closed paths $\mathcal P_1$ , $\mathcal P_2$ of 
length $2s_N$ on the complete graph on $N$ vertices $ \{1,2,\ldots, N \} $ that satisfy the following two conditions: \\
(i) $\mathcal P_1$ and $\mathcal P_2$ have at least one edge in common; \\
(ii) each edge from the union of  $\mathcal P_1$ and $\mathcal P_2$ appears at least twice in the union.\\
Indeed, if $\mathcal P_1$ and $\mathcal P_2$ do not satisfy the conditions (i) and (ii) then the corresponding term 
in the expression for $ \text{Var} \text{Tr}A_N^{2s_N} $ vanishes
due to the independence of the matrix entries up from the diagonal and the fact that the matrix entries have zero mean.
Paths $\mathcal P_1$ , $\mathcal P_2$ that satisfy (i) and (ii) are called correlated paths
(see \cite{Si-So1}, \cite{Si-So2}).   

To estimate from above the contribution of the pairs of correlated paths, 
we construct for each such pair a new path of length $4s_N-2.$   Such a mapping from the set of the pairs of 
correlated paths of length $2s_N$ to the set of paths of length $4s_N-2$ will not be injective.  In general, a path of 
length $4s_N-2$ might have many preimages.  To construct the mapping, consider an ordered pair of correlated paths
$\mathcal P_1$  and $\mathcal P_2.$  Let us consider the first edge along $\mathcal P_1$ which also belongs to 
$\mathcal P_2.$  We shall call such an edge the joint edge of the ordered pair of correlated paths 
$\mathcal P_1$  and $\mathcal P_2.$  We are now ready to construct the corresponding path of length $4s_N-2$ which will be 
denoted by $\mathcal P_1\vee \mathcal P_2.$  We choose the starting point of $\mathcal P_1\vee \mathcal P_2 $ to coincide
with the starting point of the path $\mathcal P_1.$  We begin walking along the first path until we reach for the first time
the joint edge.  At the left point of 
the joint edge we then switch to the second path.  If the directions of the joint edge in $\mathcal P_1$ and 
$\mathcal P_2$ are opposite to each other, we walk along $\mathcal P_2$ in the direction of $\mathcal P_2.$  If the 
directions of the joint edge in $\mathcal P_1$ and 
$\mathcal P_2$ coincide, we walk along $\mathcal P_2$ in the opposite direction to $\mathcal P_2.$   In both cases, we make
$2s_N-1$ steps along  $\mathcal P_2.$  In other words, we pass all $2s_N$ edges of $\mathcal P_2$ except for the joint 
edge and arrive at the right point of the joint edge.  There, we switch back to the first path and finish it. It follows 
from the construction that the new path $\mathcal P_1\vee \mathcal P_2$ is closed since it starts and ends at the starting 
point of $\mathcal P_1.$  Moreover, the length of $\mathcal P_1\vee \mathcal P_2$ is $4s_N-2$ as we omit twice the joint 
edge during our construction of $\mathcal P_1\vee \mathcal P_2.$  
We now estimate the contribution of correlated pairs $\mathcal{P}_1$, 
$\mathcal{P}_2$ in terms of $\mathcal{P}_1\vee \mathcal{P}_2$. Note that $\mathcal{P}_1 \cup \mathcal{P}_2$ and 
$\mathcal{P}_1\vee \mathcal{P}_2$ have the same edges 
appearing for the same number of times save for one important exception.
It follows from the construction of  $\mathcal{P}_1\vee \mathcal{P}_2$ that the number of appearances of the joint edge
in   $\mathcal{P}_1 \cup \mathcal{P}_2$ is bigger than the number of appearances of the joint edge in  
$\mathcal{P}_1\vee \mathcal{P}_2$  by two (in particular, if the joint edge appears only once in both $\mathcal{P}_1$
and   $\mathcal{P}_2$, it does not appear at all in $\mathcal{P}_1\vee \mathcal{P}_2$).  This observation will help us
to determine 
the number of pre-images $\mathcal{P}_1$, $\mathcal{P}_2$ of a given path $\mathcal{P}_1 \cup \mathcal{P}_2$ and relate 
the expectations associated to $\mathcal{P}_1 \cup \mathcal{P}_2$ and $\mathcal{P}_1\vee \mathcal{P}_2$.

Assume first that $\mathcal{P}_1\vee \mathcal{P}_2$ is an even path.   In this case, the arguments are identical to 
the ones used in  \cite{Si-So2} and \cite{So1}.  For the convenience of the reader, we discuss below the key steps.
To reconstruct  $\mathcal{P}_1$ and $ \mathcal{P}_2$ 
from $\mathcal{P}_1\vee \mathcal{P}_2$, it is enough to determine three things: (i) the moment of time $t_s$ in
$\mathcal{P}_1\vee \mathcal{P}_2$ where one switches from $\mathcal{P}_1$ to $\mathcal{P}_2$,  (ii)
the direction in which 
$\mathcal{P}_2$ is read, and  (iii) the origin of $\mathcal{P}_2.$  The reader can 
note that the joint edge is uniquely determined by the 
instant $t_s$, since the two endpoints of the joint edge are respectively given by the vertices occurring in 
$\mathcal{P}_1\vee \mathcal{P}_2$ at the moments $t_s$ and $t_s+2s_N-1$.
It was proved in \cite{Si-So2} (see Proposition 3) that the typical number of moments $t_s$ of possible switch is 
{\bf of} the order of $\sqrt{s_N}$ (and not $s_N$). This follows from the fact that the random walk trajectory
associated to $\mathcal{P}_1\vee \mathcal{P}_2$ does not descend below the level $x(t_s)$ during a time interval of 
length at least $2s_N.$ Given $t_s$, there are at most $2 \times 2\*s_N=4s_N$ possible choices for the orientation and origin 
of $\mathcal{P}_2.$
From that, we deduce that the contribution of correlated pairs $\mathcal{P}_1$, $ \mathcal{P}_2$ for which 
$\mathcal{P}_1\vee \mathcal{P}_2$ is an even path is of the order of
 {\bf $$ s_N^{3/2} \* \frac{1}{N} \* \frac{N}{\pi^{1/2}\* (2s_N-1)^{3/2}} \*(2\sigma)^{4s_N-2}= O((2\sigma)^{4s_N}),$$}
where the extra factor $1/N$ arises from the contribution of the erased joint edge.  Clearly, this bound 
is negligible compared to the r.h.s. of (\ref{pinsk}).

\paragraph{}We now consider the contribution of correlated paths $\mathcal P_1$ , $\mathcal P_2$ such that 
$\mathcal P_1\vee \mathcal P_2$ contains odd edges. 
To do so, we use the gluing procedure defined in \cite{Si-So2}. 
Two cases can be encountered: \begin{enumerate}
\item the joint edge of $\mathcal P_1 $ and  $\mathcal P_2$ appears in $\mathcal P_1\vee \mathcal P_2 $ exactly once
(i.e. it appears in the union of $\mathcal P_1 $ and $\mathcal P_2$ exactly three times).
\item all the odd edges of $\mathcal P_1\vee \mathcal P_2$ are read at least three times.
                              \end{enumerate}
In case 2, one can use the results established in \cite{PecheSos} to estimate the contribution to 
$\mathbb{E}[\text{Tr}M_N^{4s_N-2}]$ of paths $\mathcal{P}_1\vee \mathcal{P}_2$ admitting odd edges, all of which being 
read at least $3$ times.  Therein it is proved by using the same combinatorial machinery that proved (\ref{ocenka101})  
that $$\frac{1}{N^{2s_N-1}}\sum^*_{\mathcal P} \mathbb{E}\prod_{(m_k m_{k+1})\in \mathcal P}
|a_{m_k m_{i_{k+1}}}|\leq \frac{N}{\pi^{1/2}\*s_N^{3/2}}\*(2\sigma)^{4s_N-2}\* (1+ o(1)),$$
where the starred sum is over the set of paths $ \mathcal P $
of length $4s_N-2$ such that all odd edges (if any) are read in $ \mathcal P $ at least three times.
We first note that the number of preimages of the path $\mathcal P_1\vee \mathcal P_2$ under the described mapping 
is at most $8s_N^2.$  Indeed,  to reconstruct the pair $\mathcal P_1,\mathcal P_2$, we first note that there are
at most $2s_N$ choices for the left vertex of the joint edge of $\mathcal P_1 $ and  $\mathcal P_2$ as we select it 
among the vertices of $\mathcal P_1\vee \mathcal P_2.$  Once the left vertex of the joint edge is chosen, we recover 
the right vertex of the 
joint edge automatically since all we have to do is to make $2s_N-1$ steps along $\mathcal P_1\vee \mathcal P_2 $ 
to arrive at the right vertex of the joint edge.
Once this is done, we completely recover $\mathcal P_1. $ To recover  $\mathcal P_2,$ we  have to choose the starting 
vertex of $\mathcal P_2$ and its orientation.  This can be done in at most $2s_N \times 2= 4s_N $ ways.  Thus, we end up
with the upper bound
\begin{equation}
\label{hatyn'}
\frac{8s_N^2}{N} \* \frac{1}{N^{2s_N-1}}\sum_{\mathcal P} \mathbb{E}\prod_{(m_k m_{k+1})\in \mathcal P}
|a_{m_k m_{i_{k+1}}}|,
\end{equation}
where the sum is over the set of paths $ \mathcal P $ (i.e.  $ \mathcal P = \mathcal P_1\vee \mathcal P_2 $)
of length $4s_N-2$ such that all odd edges are read in $ \mathcal P $ at least three times (i.e. $ \mathcal P $ does not 
contain edges that appear only once there).  Using the results of \cite{PecheSos}, we can bound (\ref{hatyn'}) from above by
$$\frac{8s_N^2}{N} \* \frac{N}{\pi^{1/2}\*s_N^{3/2}}\*(2\sigma)^{4s_N-2}\* (1+ o(1)) \leq const \* \sqrt{s_N}\*
(2\sigma)^{4s_N} .$$

Finally, we have to deal with the case 1 (i.e. when the joint edge of $\mathcal P_1 $ and  $\mathcal P_2$ appears in 
$\mathcal P_1\vee \mathcal P_2 $ exactly once).
Thus we need to be able to estimate the contribution:
$$\sum^{**}_{\mathcal P} \mathbb{E}\prod_{(m_k m_{k+1})\in \mathcal P}
|a_{m_k m_{i_{k+1}}}|$$
where the two-starred sum is over the set of paths $ \mathcal P $
of length $4s_N-2$ such that all odd edges but one are read in $ \mathcal P $ at least three times.
For this, we need to modify the arguments in \cite{PecheSos} 
to include the case when there is one single edge in the path. 
We refer the reader to the above paper for the notations we will use.  As we have already pointed out,
in the case 1  the path $\mathcal P_1\vee \mathcal P_2  $ has one single edge $(ij)$, which determines  two vertices 
of the path $\mathcal P_2$. This edge serves as the joint edge of $\mathcal P_1 $ and $ \mathcal P_2.$  We recall from the 
construction of $\mathcal P_1\vee \mathcal P_2 $ that in this case, the joint edge appears three times in the union of
$\mathcal P_1 $ and $ \mathcal P_2.$  In other words, it either appears twice in $\mathcal P_1 $ and once in  
$ \mathcal P_2,$
or it appears once in $\mathcal P_1 $ and twice in  $ \mathcal P_2.$  Without loss of generality, we can assume that
the joint edge appears once in $\mathcal P_1 $ and twice in  $ \mathcal P_2.$   
Let us recall that in order to construct
 $ \mathcal P_1\vee \mathcal P_2   $, we first go along $ \mathcal P_1,$ then switch to $\mathcal P_2 $ at the appropriate 
endpoint of the joint edge, then make $2 s_N -1 $ steps along $\mathcal P_2, $ and, finally, switch back 
to $\mathcal P_1$ at the other
endpoint of the joint edge.  Let the moment of the switch back to $\mathcal P_1$  occur at time $t$ in 
$\mathcal P_1\vee \mathcal P_2 $.
Call 
$\mathcal{P}_3$ the path obtained from  $ \mathcal P_1\vee \mathcal P_2   $   by 
adding at time $t$ two successive occurrences of the (unordered) edge $(ij)$ in such a way that $\mP_3$ is still a path. 
Note that $\mP_3$  constructed in such way is a path of length $4s_N$.
Furthermore, it follows from the construction of  $\mP_3$ and the definition of the joint edge  that 
the last occurrence of $(ij)$ in $\mP_3$ is an odd edge and it necessarily starts a subsequence of odd edges  
(we refer the reader to the beginning of Section 2.1 in \cite{PecheSos} for the full account of how we split 
the set of the odd edges into disjoint subsequences $S_i, i=1, \ldots, J $ of odd edges.)
Assume that we are given a path $\mP_3$ with at least one edge read three times and where the  last two occurrences of this 
edge take place in succession. 
The idea used in \cite{PecheSos} is that a path $\mathcal{P}_3$ with odd edges (seen at least $3$ times) can be built 
from a path (or a succession of paths) with even edges  by inserting at some moments the last occurrence of odd edges. 
Given a path with even edges only, we first choose,  as described in the Insertion Procedure in Sections 3 and 4 in
\cite{PecheSos}, the set of edges that will be odd in $\mathcal{P}_3$ and choose for each of them 
the moment of time where they are going to be inserted.  To be more precise, we first recall that the set of odd edges 
can be viewed as a union of cycles.  We then split these cycles
into disjoint subsequences of odd edges to be inserted into the even path (or, in general, in the succession of paths). 
In \cite{PecheSos}, we used the 
(rough) estimate that there are at most $s_N$ possible choices for the moment of insertion of each subsequence of 
odd edges.  The expectation corresponding to such a path can be examined as in \cite{PecheSos}, 
up to the following modification. One of the subsequences $S_k$ of odd edges described in Section 2.1 of \cite{PecheSos} 
begins with the joint edge $(ij),$ and there are just two possible choices (instead of 
$s_N$)  where one can insert that particular sequence of odd edges since the moment of the insertion must follow
the moment of the appearance of $(ij).$   This follows from the fact that the edge $(ij)$ appears exactly three times 
in the path $\mP_3,$ and the last two appearances are successive.
As in \cite{PecheSos}, let us denote the number of the odd edges of $\mP_3$ by $2l.$  Since $(ij)$ is an odd edge of the path
$\mP_3,$ there are at most $2l$ ways to choose the edge $(ij)$ from the odd edges of $\mP_3.$
Once $(ij)$ is chosen, the number of the preimages of $(\mP_1, \mP_2)$
is at most $4s_N$.  Indeed, we need at most $2s_N$ choices to select the starting vertex of  $\mP_2$ and at most two 
choices to 
select the orientation of $\mP_2$. Combining these remarks, we obtain that the computations preceding Subsection 4.1.2  of 
\cite{PecheSos}
yield that 
\begin{eqnarray}
\label{baranovichi}
&&\frac{1}{N^{2s_N}}\sum^{'}_{\mathcal P_1,\mathcal P_2}\Big [\mathbb{E}
\prod_{(i_ki_{k+1})\in \mathcal P_1}\prod_{(j_kj_{k+1})\in \mathcal P'_2} |a_{i_ki_{i_{k+1}}}a_{j_kj_{k+1}}| \Big]\\
&&  \leq \frac{1}{N^{2s_N}}
\sum_{l} \frac{2}{s_N} \* (2l) \* (4s_N) \sum_{\mP_3 \text{ with $2l$ odd edges}}\mbE\prod_{(m_k m_{k+1})\in \mP_3} 
|a_{m_k m_{k+1}}|, \nonumber
\end{eqnarray}
where the sum in (\ref{baranovichi}) is over the pairs of correlated paths such that the case (1) takes place.
To apply the estimate (\ref{baranovichi}), we have to obtain an upper bound on the typical number of odd edges
in (\ref{baranovichi}).  
Thus, we only need to slightly refine our estimates given in the Subsection 4.1.2 of \cite{PecheSos}. 
As the edge $(ij)$ appears three times in $\mP_3$ and only once in $\mathcal P_1\vee \mathcal P_2,$ 
the weight of $\mP_3$ is of the order $1/N$ of the weight of 
$\mathcal P_1\vee \mathcal P_2  $ (since each matrix entry is of the order of $N^{-1/2})$.   Consider the path 
$\mathcal P_1\vee \mathcal P_2$. 
Let $\nu_N$ be the maximal number of times a vertex occurs in the even path associated to 
$\mathcal P_1\vee \mathcal P_2.$
In particular, if we know one of the endpoints of an edge (say, the left one), the number of all possible choices for the
the other endpoint is bounded from above by $\nu_N. $
Then  the number of preimages of $\mathcal P_1\vee \mathcal P_2$ is at most $\nu_N \times 4s_N.$
Indeed, since $(ij)$ is the only single edge of 
$\mathcal P_1\vee \mathcal P_2$ (i.e. the only edge appearing in
$\mathcal P_1\vee \mathcal P_2$ just once), there is no ambiguity in determining the joint edge $(ij)$ in 
$\mathcal P_1\vee \mathcal P_2.$  Then, there are at most $\nu_N$ choices to determine the place of the erased edge since we
have to select one of the appearances of the vertex $i$ in $\mathcal P_1\vee \mathcal P_2$ which can be done in at most
$\nu_N $ ways. Finally, there are $2s_N$ 
choices for the starting vertex of $\mathcal P_2$ and $2$ choices for its orientation.  As in \cite{PecheSos}, let us denote
by $\mP'$ the even path obtained from $\mathcal P_1\vee \mathcal P_2$ by the gluing procedure. 
The only modification from Subsection 4.1.2 of  \cite{PecheSos} is that the upper bound (39) on the number of ways to 
determine the cycles in the Insertion procedure  has to be 
multiplied by the factor $\nu_N^2/s_N$.  The reason for this modification is the following.  In Section 4.1.2, we observed 
that the set of odd edges can be viewed as a union of cycles.  In \cite{PecheSos}, these cycles repeat some edges of
$\mP'$.   We need to reconstruct the cycles in order to determine the set of odd edges.
Note that to reconstruct a cycle we need to know only every other 
edge  in the cycle.  For example, if we know the first and the third edges of the cycle, this uniquely determines the 
second edge of the cycle as the left endpoint of the second edge 
coincides with the right endpoint of the first edge and the right endpoint of the second edge
coincides with the left endpoint of the third edge, and so on.
In \cite{PecheSos}, we used a trivial upper bound $2s_N$ on the number of ways to choose an edge in the 
cycle since each 
such edge appears in  $\mP'$ and we have to choose it among the edges of  $\mP'.$  The difference with our situation is that 
one of 
the edges of $\mathcal P_1\vee \mathcal P_2,$ namely the joint edge $(ij),$ does not appear in $\mP'.$
However, its end points
$i$ and $j$ appear in $\mP'$ among its vertices.  Therefore, we have at most $\nu_N^2$ choices for such edge instead of 
the standard bound $2s_N$ that we used in \cite{PecheSos}.
Once the cycles are determined, one can split these cycles into disjoint sequences of odd edges to be inserted in
$\mP'$.  The total number of possible ways to insert these sequences is unchanged from Subsection 4.1.2 of  \cite{PecheSos}.
These considerations immediately imply that the contribution to the variance from the correlated paths $\mP_1, \mP_2$ is at 
most of the order 
$$ \frac{1}{N} \* \nu_N \* 4s_N \* \frac{\nu_N^2}{s_N} \* \frac{N}{s_N^{3/2}}=O(\sqrt{s_N}),$$
as long as $\nu_N <Cs_N^{2/3}.$
The case where $\nu_N>C s_N^{2/3}$ gives negligible contribution as it is extremely unlikely for any given vertex to appear 
many times in the path.  We refer the reader to Section 4.1.2 of \cite{PecheSos} where this case was analyzed.

Finally, one has to bound
$ \ \frac{1}{N^{2s_N}}\sum^{'}_{\mathcal P_1,\mathcal P_2}\Big [\mathbb{E}\prod_{(i_ki_{k+1})\in \mathcal P_1}
a_{i_ki_{i_{k+1}}}\mathbb{E}\prod_{(j_kj_{k+1})\in \mathcal P_2} a_{j_kj_{k+1}}\Big ],  \ $
where the sum is over the correlated pairs of paths.  This can be done in the same way as we treated
$ \ \frac{1}{N^{2s_N}}\sum^{'}_{\mathcal P_1,\mathcal P_2}\Big [\mathbb{E}
\prod_{(i_ki_{k+1})\in \mathcal P_1}\prod_{(j_kj_{k+1})\in \mathcal P'_2} a_{i_ki_{i_{k+1}}}a_{j_kj_{k+1}}\Big] \ $ above.
This finishes the proof of the lemma and {\bf gives us} the proof of the main result.

\end{document}